\documentclass[12pt]{article}

\usepackage{amsmath,amssymb,amsthm}

\begin{document}

\newcommand{\End}{{\rm{End}\ts}}
\newcommand{\non}{\nonumber}
\newcommand{\wt}{\widetilde}
\newcommand{\wh}{\widehat}
\newcommand{\ot}{\otimes}
\newcommand{\la}{\lambda}
\newcommand{\al}{\alpha}
\newcommand{\be}{\beta}
\newcommand{\ga}{\gamma}
\newcommand{\si}{\sigma}
\newcommand{\vs}{\varsigma}
\newcommand{\de}{\delta^{}}
\newcommand{\De}{\Delta}
\newcommand{\Om}{\Omega}
\newcommand{\ze}{\zeta}
\newcommand{\om}{\omega^{}}
\newcommand{\hra}{\hookrightarrow}
\newcommand{\ve}{\varepsilon}
\newcommand{\ts}{\,}
\newcommand{\qin}{q^{-1}}
\newcommand{\tss}{\hspace{1pt}}
\newcommand{\U}{ {\rm U}}
\newcommand{\Y}{ {\rm Y}}
\newcommand{\CC}{\mathbb{C}\tss}
\newcommand{\RR}{\mathbb{R}\tss}
\newcommand{\ZZ}{\mathbb{Z}}
\newcommand{\A}{\mathcal{A}}
\newcommand{\Z}{{\rm Z}}
\newcommand{\gl}{\mathfrak{gl}}
\newcommand{\oa}{\mathfrak{o}}
\newcommand{\spa}{\mathfrak{sp}}
\newcommand{\g}{\mathfrak{g}}
\newcommand{\ka}{\mathfrak{k}}
\newcommand{\p}{\mathfrak{p}}
\newcommand{\sll}{\mathfrak{sl}}
\newcommand{\agot}{\mathfrak{a}}
\newcommand{\qdet}{ {\rm qdet}\ts}
\newcommand{\sdet}{ {\rm sdet}\ts}
\newcommand{\sgn}{ {\rm sgn}}
\newcommand{\Sym}{\mathfrak S}
\newcommand{\Fand}{\qquad\text{and}\qquad}

\renewcommand{\theequation}{\arabic{section}.\arabic{equation}}

\newtheorem{thm}{Theorem}[section]
\newtheorem{lem}[thm]{Lemma}
\newtheorem{prop}[thm]{Proposition}
\newtheorem{cor}[thm]{Corollary}

\theoremstyle{definition}
\newtheorem{defin}[thm]{Definition}
\newtheorem{example}[thm]{Example}

\theoremstyle{remark}
\newtheorem{remark}[thm]{Remark}

\newcommand{\bth}{\begin{thm}}
\renewcommand{\eth}{\end{thm}}
\newcommand{\bpr}{\begin{prop}}
\newcommand{\epr}{\end{prop}}
\newcommand{\ble}{\begin{lem}}
\newcommand{\ele}{\end{lem}}
\newcommand{\bco}{\begin{cor}}
\newcommand{\eco}{\end{cor}}
\newcommand{\bde}{\begin{defin}}
\newcommand{\ede}{\end{defin}}
\newcommand{\bex}{\begin{example}}
\newcommand{\eex}{\end{example}}
\newcommand{\bre}{\begin{remark}}
\newcommand{\ere}{\end{remark}}

\newcommand{\bal}{\begin{aligned}}
\newcommand{\eal}{\end{aligned}}
\newcommand{\beq}{\begin{equation}}
\newcommand{\eeq}{\end{equation}}
\newcommand{\ben}{\begin{equation*}}
\newcommand{\een}{\end{equation*}}

\newcommand{\bpf}{\begin{proof}}
\newcommand{\epf}{\end{proof}}

\def\beql#1{\begin{equation}\label{#1}}

\title{\Large\bf  Representations of the twisted quantized enveloping
algebra of type $C_n$}
\author{{\sc A. I. Molev}\\[15pt]
School of Mathematics and Statistics\\
University of Sydney,
NSW 2006, Australia\\
{\tt alexm\hspace{0.09em}@\hspace{0.1em}maths.usyd.edu.au}}

\date{} 

\maketitle

\vspace{12 mm}

\begin{abstract}
We prove a version
of the Poincar\'e--Birkhoff--Witt theorem
for the twisted quantized enveloping algebra
$\U'_q(\spa_{2n})$.
This is a subalgebra of
$\U_q(\gl_{2n})$
and a deformation of the universal enveloping algebra
$\U(\spa_{2n})$ of the symplectic Lie algebra.
We classify finite-dimensional irreducible
representations of
$\U'_q(\spa_{2n})$ in terms of their highest weights
and show that these representations are deformations
of the finite-dimensional irreducible
representations of $\spa_{2n}$.

\vspace{7 mm}

{\it Mathematics Subject Classification: {\rm 81R10}}

\end{abstract}

\newpage

\section{Introduction}\label{sec:int}
\setcounter{equation}{0}

Let $\g$ denote the orthogonal or symplectic Lie algebra
$\oa_N$ or $\spa_{2n}$ over the field of complex numbers.
There are at least two different $q$-analogues of the
universal enveloping algebra $\U(\g)$. These are
the {\it quantized enveloping algebra\/} $\U_q(\g)$
introduced by Drinfeld~\cite{d:ha} and
Jimbo~\cite{j:qd}, and
the {\it twisted\/} (or {\it nonstandard\/})
{\it quantized enveloping algebra\/} $\U'_q(\g)$,
introduced by Gavrilik and Klimyk~\cite{gk:qd} for $\g=\oa_N$
and by Noumi~\cite{n:ms} for $\g=\spa_{2n}$.
If $q$ is a complex number which is nonzero and not a root of unity
then the representation theory of $\U_q(\g)$ is very much similar
to that of the Lie algebra $\g$; see e.g.
Chari and Pressley~\cite[Chapter~10]{cp:gq}.
The description of the finite-dimensional irreducible representations
of $\U'_q(\oa_N)$ given by Iorgov and Klimyk~\cite{ik:ct}
both exhibits similarity with the classical theory and reveals
some new phenomena specific to the quantum case.

In this paper we are concerned with the description
of finite-dimensional irreducible representations
of the twisted quantized enveloping algebra $\U'_q(\spa_{2n})$.
We introduce a class of
highest weight representations for this algebra and show
that any finite-dimensional irreducible representation
of $\U'_q(\spa_{2n})$ is highest weight. Then we give
necessary and sufficient conditions
on a highest weight representation
to be finite-dimensional. By the main theorem (Theorem~\ref{thm:classif}),
the finite-dimensional irreducible representations
of $\U'_q(\spa_{2n})$ are naturally parameterized
by the $n$-tuples
\ben
\la=(\si_1\ts q^{m_1},\dots,\si_n\ts q^{m_n}),
\een
where the $m_i$ are positive integers satisfying
$m_1\leqslant m_2\leqslant\dots\leqslant m_n$ and each
$\si_i$ equals $1$ or $-1$.
Regarding this parametrization, the algebra $\U'_q(\spa_{2n})$
appears to be closer to
the quantized enveloping algebra $\U_q(\spa_{2n})$ than
to its orthogonal counterpart $\U'_q(\oa_N)$.

We work with the presentation of $\U'_q(\spa_{2n})$ introduced
in \cite{mrs:cs} which is a slight modification of
Noumi's definition \cite{n:ms}.
The defining relations are written in a matrix form
as a reflection equation for the matrix of generators.
A version of the Poincar\'e--Birkhoff--Witt theorem
for this algebra over the field $\CC(q)$ was proved
in \cite{mrs:cs}. Here we consider $q$ to be
a nonzero complex number such that $q^2\ne 1$ and prove
the corresponding theorem for $\U'_q(\spa_{2n})$ (Theorem~\ref{thm:pbw}),
relying on the Poincar\'e--Birkhoff--Witt theorem
for the quantized enveloping algebra $\U_q(\gl_{N})$.

Note that a similar reflection equation presentation exists for the
algebra $\U'_q(\oa_N)$ as well. These presentations
were derived in \cite{n:ms}
by regarding $\U'_q(\oa_N)$ and $\U'_q(\spa_{2n})$
as subalgebras of $\U_q(\gl_N)$ for appropriate $N$;
see also Gavrilik, Iorgov and Klimyk~\cite{gik:nd}.
These subalgebras are coideals with respect
to the coproduct on $\U_q(\gl_N)$.
A more general description of the coideal subalgebras of
the quantized enveloping algebras
was given by Letzter~\cite{l:sp, l:cs}.
We outline a new proof of the
Poincar\'e--Birkhoff--Witt theorem
for $\U'_q(\oa_N)$ analogous to the symplectic case
(see Remark~\ref{rem:osyquo} below);
cf. Iorgov and Klimyk~\cite{ik:nd}.

The algebra $\U'_q(\oa_N)$ and its representations were studied
by many authors. In particular, it plays the role
of the symmetry
algebra for the $q$-oscillator representation of
the quantized enveloping algebra $\U_q(\spa_{2n})$;
see Noumi, Umeda and Wakayama~\cite{nuw:dp}.
A quantum analogue of the Brauer algebra associated with
$\U'_q(\oa_N)$ was constructed in \cite{m:nq}.
Some families of Casimir elements
were produced by Noumi, Umeda and Wakayama~\cite{nuw:dp},
Havl\'\i\v cek, Klimyk and Po\v sta~\cite{hkp:ce}
and by Gavrilik and Iorgov~\cite{gi:ce} for the algebra $\U'_q(\oa_N)$,
and by Molev, Ragoucy and Sorba~\cite{mrs:cs} for both
$\U'_q(\oa_N)$ and $\U'_q(\spa_{2n})$. The paper \cite{mrs:cs} also provides
a construction of certain Yangian-type algebras associated with
the twisted quantized enveloping algebras.
Their applications to spin chain models were discussed
in Arnaudon {\it et al.}~\cite{acdfr:sb}.
The algebra $\U'_q(\spa_{2n})$ appears to
have received much less attention in the literature as compared
to its orthogonal counterpart, which we hope to remedy by this work.

\section{Definitions and preliminaries}\label{sec:dp}
\setcounter{equation}{0}

Fix a complex parameter $q$ which is nonzero and $q^2\ne1$.
Following Jimbo~\cite{j:qu},
we introduce the $q$-analogue $\U_q(\gl_N)$ of the universal enveloping
algebra $\U(\gl_N)$ as an associative algebra generated by
the elements $t_1,\dots,t_N,t_1^{-1},\dots,t_N^{-1}$, $e_1,\dots,e_{N-1}$ and
$f_1,\dots,f_{N-1}$ with the defining relations
\ben
\bal
t_i\tss t_j=t_j\tss t_i, &\qquad t^{}_i\tss t_i^{-1}=t_i^{-1}\tss t^{}_i=1, \\
t^{}_i\tss e^{}_j\tss t_i^{-1}=e^{}_j \ts q^{\tss\delta_{ij}-\delta_{i,j+1}},&\qquad
t^{}_i\tss f^{}_j\tss t_i^{-1}=f^{}_j \ts q^{-\delta_{ij}+\delta_{i,j+1}},\\
[e_i,f_j]=\delta_{ij}\ts\frac{k^{}_i-k_i^{-1}}{q-\qin}&\qquad \text{with}\ \
k_i=t^{}_i\tss t_{i+1}^{-1},\\
[e_i,e_j]=[f_i,f_j]=0&\qquad\text{if}\ \ |i-j|>1,\\
e^2_i\tss e^{}_{j}-(q+\qin)\tss e^{}_i\tss e^{}_{j}\tss e^{}_i
&+e^{}_{j}\tss e^2_i=0\qquad\text{if}\ \ |i-j|=1,\\
f^2_i\tss f^{}_{j}-(q+\qin)\tss f^{}_i\tss f^{}_{j}\tss f^{}_i
&+f^{}_{j}\tss f^2_i=0\qquad\text{if}\ \ |i-j|=1.
\eal
\een
The $q$-analogues of the root vectors are defined inductively by
\begin{alignat}{2}
{}&e_{i,i+1}=e_i,\qquad {}&&e_{i+1,i}=f_i,
\non\\
\label{rootv}
{}&e_{ij}=e_{ik}\tss e_{kj}-q\ts e_{kj}\tss e_{ik}\qquad &&\text{for}\ \ i<k<j,\\
{}&e_{ij}=e_{ik}\tss e_{kj}-\qin\tss e_{kj}\tss e_{ik}\qquad &&\text{for}\ \ i>k>j,
\non
\end{alignat}
and the elements $e_{ij}$ are independent of the choice of values
of the index $k$.

Following \cite{j:qu} and \cite{rtf:ql}, consider
the $R$-matrix presentation of the algebra $\U_q(\gl_N)$.
The $R$-matrix is given by
\beql{rmatrixc}
R=q\ts\sum_i E_{ii}\ot E_{ii}+\sum_{i\ne j} E_{ii}\ot E_{jj}+
(q-\qin)\sum_{i<j}E_{ij}\ot E_{ji}
\eeq
which is an element of $\End\CC^N\ot \End\CC^N$, where
the $E_{ij}$ denote the standard matrix units and the indices run over
the set $\{1,\dots,N\}$. The $R$-matrix satisfies the Yang--Baxter equation
\beql{YBEconst}
R_{12}\ts R_{13}\ts  R_{23} =  R_{23}\ts  R_{13}\ts  R_{12},
\eeq
where both sides take values in $\End\CC^N\ot \End\CC^N\ot \End\CC^N$ and
the subindices indicate the copies of $\End\CC^N$, e.g.,
$R_{12}=R\ot 1$ etc.

The quantized enveloping algebra $\U_q(\gl_N)$ is generated
by elements $t_{ij}$ and $\bar t_{ij}$ with $1\leqslant i,j\leqslant N$
subject to the relations
\beql{defrel}
\bal
t_{ij}&=\bar t_{ji}=0, \qquad 1 \leqslant i<j\leqslant N,\\
t_{ii}\ts \bar t_{ii}&=\bar t_{ii}\ts t_{ii}=1,\qquad 1\leqslant i\leqslant N,\\
R\ts T_1T_2&=T_2T_1R,\qquad R\ts \overline T_1\overline T_2=
\overline T_2\overline T_1R,\qquad
R\ts \overline T_1T_2=T_2\overline T_1R.
\eal
\eeq
Here $T$ and $\overline T$ are the matrices
\beql{matrt}
T=\sum_{i,j}t_{ij}\ot E_{ij},\qquad \overline T=\sum_{i,j}
\overline t_{ij}\ot E_{ij},
\eeq
which are regarded as elements of the algebra $\U_q(\gl_N)\ot \End\CC^N$.
Both sides of each of the $R$-matrix relations in \eqref{defrel}
are elements of $\U_q(\gl_N)\ot \End\CC^N\ot \End\CC^N$ and the subindices
of $T$ and $\overline T$ indicate the copies of $\End\CC^N$ where
$T$ or $\overline T$ acts; e.g. $T_1=T\ot 1$. In terms of the
generators the defining relations between the $t_{ij}$
can be written as
\beql{defrelg}
q^{\delta_{ij}}\ts t_{ia}\ts t_{jb}-
q^{\delta_{ab}}\ts t_{jb}\ts t_{ia}
=(q-\qin)\ts (\de_{b<a} -\de_{i<j})
\ts t_{ja}\ts t_{ib},
\eeq
where $\de_{i<j}$ equals $1$ if $i<j$ and $0$ otherwise.
The relations between the $\bar t_{ij}$
are obtained by replacing $t_{ij}$ by $\bar t_{ij}$ everywhere in
\eqref{defrelg} and the relations involving both
$t_{ij}$ and $\bar t_{ij}$ have the form
\beql{defrelg2}
q^{\delta_{ij}}\ts \bar t_{ia}\ts t_{jb}-
q^{\delta_{ab}}\ts t_{jb}\ts \bar t_{ia}
=(q-\qin)\ts (\de_{b<a}\ts t_{ja}\ts \bar t_{ib} -\de_{i<j}\ts
\ts \bar t_{ja}\ts t_{ib}).
\eeq

An isomorphism between the two presentations is given by the formulas
\beql{isomgln}
t_i\mapsto t_{ii},\qquad t^{-1}_i\mapsto \bar t_{ii},\qquad
e_{ij}\mapsto -\frac{\bar t_{ij}\ts t_{ii}}{q-\qin},\qquad
e_{ji}\mapsto \frac{\bar t_{ii}\ts t_{ji}}{q-\qin},
\end{equation}
for $i<j$; see e.g. \cite[Section~8.5.5]{ks:qg}.
We shall identify the corresponding elements of $\U_q(\gl_N)$
via this isomorphism.

For any $N$-tuple $\vs=(\vs_1,\dots,\vs_N)$
with $\vs_i\in\{-1,1\}$ the mapping
\beql{autosi}
t_{ij}\mapsto \vs_i\ts t_{ij},\qquad
\bar t_{ij}\mapsto \vs_i\ts \bar t_{ij},
\eeq
defines an automorphism of the algebra $\U_q(\gl_N)$.

By the Poincar\'e--Birkhoff--Witt theorem for the
$A_n$ type quantized enveloping algebra,
the monomials
\begin{multline}
e^{k_{N,N-1}}_{N,N-1}\ts e^{k_{N,N-2}}_{N,N-2}\ts
e^{k_{N-1,N-2}}_{N-1,N-2}
\dots\ts e^{k_{N2}}_{N2}
\dots\ts e^{k_{32}}_{32}\ts e^{k_{N1}}_{N1}\ts
\dots\ts e^{k_{21}}_{21}\\
{}\times t_1^{\ts l_1}\dots\ts t_N^{\ts l_N}\ts
e^{k_{12}}_{12}\dots\ts e^{k_{1N}}_{1N}\ts
e^{k_{23}}_{23}\dots\ts e^{k_{2N}}_{2N}\dots\ts
e^{k_{N-1,N}}_{N-1,N},
\end{multline}
where the $k_{ij}$ run over non-negative integers
and the $l_i$ run over integers,
form a basis of the algebra of $\U_q(\gl_N)$;
see \cite{r:ap}, \cite[Proposition~9.2.2]{cp:gq}
(this basis corresponds to the reduced decomposition
\ben
w_0=s^{}_{N-1}\ts s^{}_{N-2}\ts s^{}_{N-1}\dots\ts
s^{}_2\ts s^{}_3\dots\ts s^{}_{N-1}\ts s^{}_1\ts s^{}_2\dots\ts s^{}_{N-1}
\een
of the longest element of the Weyl group).
Using the isomorphism \eqref{isomgln}, we may conclude that
the monomials
\begin{multline}
t^{k_{N,N-1}}_{N,N-1}\ts t^{k_{N,N-2}}_{N,N-2}\ts
t^{k_{N-1,N-2}}_{N-1,N-2}
\dots\ts t^{k_{N2}}_{N2}
\dots\ts t^{k_{32}}_{32}\ts t^{k_{N1}}_{N1}\ts
\dots\ts t^{k_{21}}_{21}\\
{}\times t_{11}^{l_1}\dots t_{NN}^{l_N}\ts
\bar t^{\ts k_{12}}_{12}\dots\ts \bar t^{\ts k_{1N}}_{1N}\ts
\bar t^{\ts k_{23}}_{23}\dots\ts \bar t^{\ts k_{2N}}_{2N}\dots\ts
\bar t^{\ts k_{N-1,N}}_{N-1,N},
\end{multline}
where the $k_{ij}$ run over non-negative integers
and the $l_i$ run over integers, also
form a basis of $\U_q(\gl_N)$.
This is implied by the relations
\beql{tiicom}
t_{ii}\ts t_{jb}=q^{\de_{ij}-\de_{ib}}\ts t_{jb}\ts t_{ii}\Fand
t_{ii}\ts \bar t_{jb}=q^{\de_{ij}-\de_{ib}}\ts \bar t_{jb}\ts t_{ii}.
\eeq

Let $\U^-$ denote the subalgebra of $\U_q(\gl_N)$ generated
by the elements $t_{ii}$ with $i=1,\dots,N$ and
$t_{ij}$ with $1\leqslant j<i\leqslant N$.
Fix a permutation $\pi$ of the indices $2,3,\dots,N$ and
consider a linear ordering $\prec$ on the set of elements
$t_{ij}$ with $1\leqslant j<i\leqslant N$ such that $\pi(i)<\pi(k)$
implies $t_{ij}\prec t_{kl}$ for all possible $j$ and $l$.

\bpr\label{prop:arbord}
For the linear ordering defined as above,
the ordered monomials of the form
\beql{lengthmon}
\prod_{i>j}t_{ij}^{k_{ij}}\ts\prod_i t_{ii}^{m_i},
\eeq
where the $k_{ij}$ and $m_i$ are non-negative integers,
form a basis of $\U^-$.
\epr

\bpf
For any non-negative integer $l$ consider
the subspace $\U^-_l$ of $\U^-$
of elements of degree at most $l$
in the generators. Due to the Poincar\'e--Birkhoff--Witt theorem,
a basis of $\U^-_l$ is formed by the monomials
\ben
t^{k_{N,N-1}}_{N,N-1}\ts t^{k_{N,N-2}}_{N,N-2}\ts
t^{k_{N-1,N-2}}_{N-1,N-2}
\dots\ts t^{k_{N2}}_{N2}
\dots\ts t^{k_{32}}_{32}\ts t^{k_{N1}}_{N1}\ts
\dots\ts t^{k_{21}}_{21}\ts t^{m_{1}}_{11}\ts
\dots\ts t^{m_{N}}_{NN},
\een
with the sum of all powers not exceeding $l$.
Hence, it will be sufficient to show that
the ordered monomials \eqref{lengthmon}
span $\U^-$.
The statement will then follow by counting the number
of the ordered monomials of degree not exceeding $l$.

By the defining relations \eqref{defrelg}, we have
\ben
t_{ia}\ts t_{jb}=t_{jb}\ts t_{ia}
+(q-\qin)\ts t_{ja}\ts t_{ib},\qquad t_{ja}\ts t_{ib}=t_{ib}\ts t_{ja}
\een
for $i>j>a>b$, while
\beql{relco}
t_{ia}\ts t_{ab}=t_{ab}\ts t_{ia}
+(q-\qin)\ts t_{aa}\ts t_{ib},\qquad
t_{ia}\ts t_{ib}=q\ts t_{ib}\ts t_{ia}
\eeq
for $i>a>b$. Given a monomial
$t_{i_1a_1}\cdots\ts t_{i_pa_p}$ with $i_r>a_r$ for all $r$, an easy induction
on the degree $p$ shows that modulo elements of degree ${}<p$,
this monomial can be written as a linear combination
of monomials of the form $t_{j_1b_1}\cdots\ts t_{j_pb_p}$ where
$\pi(j_1)\leqslant\cdots \leqslant \pi(j_p)$.
By the second relation in \eqref{relco}, this monomial coincides
with an ordered monomial up to a factor which is a power of $q$.
The proof is completed by taking
the first relation in \eqref{tiicom} into account.
\epf

We shall also use an extended quantized enveloping algebra
$\hat\U_q(\gl_N)$.
This is an associative algebra generated
by elements $t_{ij}$ and $\bar t_{ij}$ with $1\leqslant i,j\leqslant N$
and elements $t^{-1}_{ii}$ and $\bar t_{ii}^{\ts\ts-1}$ with $1\leqslant i\leqslant N$
subject to the relations
\ben
\bal
t_{ij}&=\bar t_{ji}=0, \qquad 1 \leqslant i<j\leqslant N,\\
t_{ii}\ts \bar t_{ii}&=\bar t_{ii}\ts t_{ii},\qquad
t^{}_{ii}\ts t^{-1}_{ii}=t^{-1}_{ii}\ts t^{}_{ii}=1,\qquad
\bar t^{}_{ii}\ts \bar t_{ii}^{\ts\ts-1}=
\bar t^{\ts\ts-1}_{ii}\ts \bar t_{ii}^{}=1,
\qquad 1\leqslant i\leqslant N,\\
R\ts T_1T_2&=T_2T_1R,\qquad R\ts \overline T_1\overline T_2=
\overline T_2\overline T_1R,\qquad
R\ts \overline T_1T_2=T_2\overline T_1R,
\eal
\een
where we use the notation of \eqref{defrel}.
Obviously,
we have a natural epimorphism $\hat\U_q(\gl_N)\to\U_q(\gl_N)$
which takes the generators $t_{ij}$ and $\bar t_{ij}$ of
$\hat\U_q(\gl_N)$ respectively to the elements of $\U_q(\gl_N)$
with the same name. More generally, for any nonzero complex numbers
$\rho_i$ with $i=1,\dots,N$ the mapping
\beql{epimor}
\varrho:
t_{ij}\mapsto \rho_i\ts t_{ij},\qquad
\bar t_{ij}\mapsto \rho_i\ts \bar t_{ij}
\eeq
defines an  epimorphism $\hat\U_q(\gl_N)\to\U_q(\gl_N)$.
For any $i\in\{1,\dots,N\}$ the element $t_{ii}\tss \bar t_{ii}$
belongs to the center of $\hat\U_q(\gl_N)$
while $t_{ii}\tss \bar t_{ii}-\rho_i^2$ is contained
in the kernel of $\varrho$.

Let $\hat\U^0$ denote the (commutative) subalgebra of
$\hat\U_q(\gl_N)$ generated
by the elements $t_{ii}$, $t^{-1}_{ii}$, $\bar t_{ii}$ and
$\bar t_{ii}^{\ts\ts-1}$ with $i=1,\dots,N$,
and let $\hat\U^-$ denote the subalgebra of $\hat\U_q(\gl_N)$
generated by $\hat\U^0$ and all elements
$t_{ij}$ with $i>j$.

Fix a permutation $\pi$ of the indices $2,3,\dots,N$ and
consider a linear ordering $\prec$ on the set of elements
$t_{ij}$ with $1\leqslant j<i\leqslant N$ such that $\pi(i)<\pi(k)$
implies $t_{ij}\prec t_{kl}$ for all possible $j$ and $l$.

\bco\label{cor:extpbw}
The subalgebra $\hat\U^0$ is isomorphic
to the algebra of Laurent polynomials in
the variables $t_{ii}$ and $\bar t_{ii}$.
Moreover, for the linear ordering defined as above,
the ordered monomials in the elements
$t_{ij}$ with $1\leqslant j<i\leqslant N$,
form a basis of the right $\hat\U^0$-module $\hat\U^-$.
\eco

\bpf
Obviously, the subalgebra $\hat\U^0$ is spanned by
the Laurent monomials in the elements $t_{ii}$ and $\bar t_{ii}$.
We need to verify that the monomials are linearly independent.
Suppose that
\ben
\sum_{\mathbf m,\mathbf l} c_{\mathbf m,\mathbf l}\ts
t_{11}^{\tss m_1}\bar t_{11}^{\ts l_1}\dots\ts
t_{NN}^{\tss m_N}\bar t_{NN}^{\ts l_N}=0,
\een
summed over $N$-tuples of integers
$\mathbf m=(m_1,\dots,m_N)$ and $\mathbf l=(l_1,\dots,l_N)$, where
only a finite
number of the coefficients $c_{\mathbf m,\mathbf l}\in\CC$ is nonzero.
Apply an epimorphism of the form \eqref{epimor} to both sides
of this relation. This gives a relation in $\U^0$,
\ben
\sum_{\mathbf m,\mathbf l} c_{\mathbf m,\mathbf l}\ts
\rho_1^{\tss m_1+l_1}\dots\ts \rho_N^{\tss m_N+l_N}
t_{11}^{\tss m_1-l_1}\dots\ts
t_{NN}^{\tss m_N-l_N}=0.
\een
Since the monomials $t_{11}^{\tss r_1}\dots\ts
t_{NN}^{\tss r_N}$ with $r_i\in\ZZ$ are linearly independent,
we get
\ben
\sum_{\mathbf m,\mathbf l} c_{\mathbf m,\mathbf l}\ts
\rho_1^{\tss m_1+l_1}\dots\ts \rho_N^{\tss m_N+l_N}=0
\een
for any fixed integer differences $m_i-l_i$ with $i=1,\dots, N$.
Varying the values of the parameters $\rho_i$ we conclude that
$c_{\mathbf m,\mathbf l}=0$ for all $\mathbf m$ and $\mathbf l$.
This proves the first part of the corollary.

Arguing as in the proof of Proposition~\ref{prop:arbord},
we obtain that
the ordered monomials in the elements
$t_{ij}$ with $i>j$ span $\hat\U^-$ as a right $\hat\U^0$-module.
It remains to show that the ordered monomials are linearly
independent over $\hat\U^0$. Suppose that
a linear combination of the ordered monomials is zero,
\ben
\sum_{\mathbf k}
t_{i_1a_1}^{k_{i_1a_1}}\dots\ts
t_{i_ma_m}^{k_{i_ma_m}}\ts d_{\mathbf k}=0,\qquad d_{\mathbf k}\in\hat\U^0,
\een
where $\{(i_1a_1),\dots, (i_m,a_m)\}=\{(2,1),(3,1),(3,2),\dots, (N,N-1)\}$
and $\mathbf k$ runs over
a finite set of tuples of non-negative integers $k_{i_ra_r}$.
Apply an epimorphism of the form \eqref{epimor} to both sides
of this relation. By Proposition~\ref{prop:arbord},
we get $\varrho(d_{\mathbf k})=0$ for all $\mathbf k$.
As in the first part of the proof, varying the parameters $\rho_i$
we conclude that $d_{\mathbf k}=0$ for all $\mathbf k$.
\epf

Note that a similar argument can be used to demonstrate that
the ordered monomials of the form
\ben
t_{i_1a_1}^{k_{i_1a_1}}\dots\ts
t_{i_ma_m}^{k_{i_ma_m}}\ts \bar t_{a_mi_m}^{\ts k_{a_mi_m}}\dots\ts
\bar t_{a_1i_1}^{\ts k_{a_1i_1}}
\een
form a basis of the left or right $\hat\U^0$-module
$\hat\U_q(\gl_N)$.

Now we reproduce some results from \cite{mrs:cs} and \cite{n:ms}
concerning the twisted
quantized enveloping algebra $\U'_q(\spa_{2n})$.
This is an associative algebra generated by elements $s_{ij}$,
$i,j\in\{1,\dots,2n\}$ and $s_{i,i+1}^{-1}$, $i=1,3,\dots,2n-1$.
The generators $s_{ij}$ are zero
for $j=i+1$ with even $i$, and
for $j\geqslant i+2$ and all $i$. We combine the $s_{ij}$ into
a matrix $S$ as in \eqref{matrt},
\beql{matrts}
S=\sum_{i,j}s_{ij}\ot E_{ij},
\eeq
so that $S$ has a block-triangular form with
$n$ diagonal $2\times 2$-blocks,
\ben
S=\begin{pmatrix} s_{11}&s_{12}&0&0&\cdots&0&0\\
                     s_{21}&s_{22}&0&0&\cdots&0&0\\
                     s_{31}&s_{32}&s_{33}&s_{34}&\cdots&0&0\\
                     s_{41}&s_{42}&s_{43}&s_{44}&\cdots&0&0\\
                         \vdots&\vdots&\vdots&\vdots&\ddots&\vdots&\vdots\\
                     s_{2n-1,1}&s_{2n-1,2}
                         &s_{2n-1,3}&s_{2n-1,4}&
                                \cdots&s_{2n-1,2n-1}&s_{2n-1,2n}\\
                     s_{2n,1}&s_{2n,2}&s_{2n,3}&s_{2n,4}&
                                               \cdots&s_{2n,2n-1}&s_{2n,2n}
\end{pmatrix}.
\een
Consider the transposed
to the $R$-matrix \eqref{rmatrixc} given by
\beql{rt}
R^{\tss\prime}=q\ts\sum_i E_{ii}\ot E_{ii}+\sum_{i\ne j} E_{ii}\ot E_{jj}+
(q-\qin)\sum_{i<j}E_{ji}\ot E_{ji}.
\eeq
The defining relations of $\U'_q(\spa_{2n})$ have the form
of a reflection equation
\beql{refsymp}
R\ts S_1 R^{\tss\prime} S_2=S_2R^{\tss\prime} S_1R,
\eeq
together with
\beql{invrel}
s^{}_{i,i+1}\ts s_{i,i+1}^{-1}=s_{i,i+1}^{-1}\ts s^{}_{i,i+1}=1
\eeq
and
\beql{qdetrel}
s^{}_{i+1,i+1}\ts s^{}_{ii}-q^2\ts s^{}_{i+1,i}s^{}_{i,i+1}=q^3
\eeq
for $i=1,3,\dots,2n-1$. In terms of the generators $s_{ij}$ the relation
\eqref{refsymp} is written as
\begin{align}\label{drabss}
q^{\delta_{aj}+\delta_{ij}}\ts s_{ia}\ts s_{jb}-
q^{\delta_{ab}+\delta_{ib}}\ts s_{jb}\ts s_{ia}
{}&=(q-\qin)\ts q^{\delta_{ai}}\ts (\de_{b<a} -\de_{i<j})
\ts s_{ja}\ts s_{ib}\\
{}&+(q-\qin)\ts \big(q^{\delta_{ab}}\ts \de_{b<i}\ts s_{ji}\ts s_{ba}
- q^{\delta_{ij}}\ts \de_{a<j}\ts s_{ij}\ts s_{ab}\big)
\non\\
{}&+ (q-\qin)^2\ts  (\de_{b<a<i} -\de_{a<i<j})\ts s_{ji}\ts s_{ab},
\non
\end{align}
where $\de_{i<j}$ or $\de_{i<j<k}$ equals $1$ if the subindex inequality
is satisfied and $0$ otherwise.

For any $2n$-tuple $\vs=(\vs_1,\dots,\vs_{2n})$
with $\vs_i\in\{-1,1\}$ the mapping
\beql{autosic}
s_{ij}\mapsto \vs_i\ts\vs_j\ts s_{ij},
\eeq
defines an automorphism of the algebra $\U'_q(\spa_{2n})$.
This is verified directly from the
defining relations of the algebra.

Introduce the $2n\times 2n$ matrix $G$ by
\beql{g}
G=q\ts \sum_{k=1}^n E_{2k-1,2k}-\sum_{k=1}^n E_{2k,2k-1}
\eeq
so that
\ben
G=\left(\begin{matrix}
                     0&q&\cdots&0&0\\
                     -1&0&\cdots&0&0\\
                     \vdots&\vdots&\ddots&\vdots&\vdots\\
                     0&0&\cdots&0&q\\
                     0&0&\cdots&-1&0
\end{matrix}\right).
\een
The mapping
\beql{atemb}
S\mapsto T\ts G\ts \overline T^{\ts t}
\eeq
defines a homomorphism $\U'_q(\spa_{2n})\to\U_q(\gl_{2n})$.
Explicitly, in terms of generators it is written as
\beql{sijtij}
s_{ij}\mapsto q\ts \sum_{k=1}^n
t_{i,2k-1}\ts\bar t_{j,2k}
-\sum_{k=1}^n t_{i,2k}\ts\bar t_{j,2k-1}.
\eeq

When $q$ is regarded as a variable,
the mapping \eqref{atemb} is an embedding
of the algebras over the field $\CC(q)$;
see \cite[Theorem~2.8]{mrs:cs}.
This is implied by the fact that the algebra
$\U'_q(\spa_{2n})$ specializes to $\U(\spa_{2n})$ as $q\to 1$.
More precisely, set $\A=\CC[q,\qin]$ and consider
the $\A$-subalgebra $\U'_{\A}$ of $\U'_q(\spa_{2n})$
generated by the elements
\ben
\sigma_{i,i+1}=\frac{s_{i,i+1}-q}{q-1},\qquad
\sigma_{i+1,i}=\frac{s_{i+1,i}+1}{q-1},\qquad i=1,3,\dots,2n-1,
\een
and
\ben
\sigma_{ij}=\frac{s_{ij}}{q-\qin},\qquad i\geqslant j,
\quad\text{excluding $i=j+1$, $j$ odd}.
\een
Then we have an isomorphism
\beql{speci}
\U'_{\A}{\ot}^{}_{\A}\ts\CC\cong \U(\spa_{2n}),
\eeq
where the action of $\A$ on
$\CC$ defined via the evaluation $q=1$; see \cite{mrs:cs}.
The symplectic Lie algebra $\spa_{2n}$ is defined as
the subalgebra of $\gl_{2n}$ spanned by the elements
\ben
F_{ij}=\sum_{k=1}^{2n} (E_{ik}\ts g_{kj}+E_{jk}\ts g_{ki}),
\een
where the $g_{ij}$ are the matrix elements of the matrix
$G^{\tss\circ}=[\tss g_{ij}]$
obtained by evaluating $G$ at $q=1$,
\ben
G^{\tss\circ}=\left(\begin{matrix}
                     0&1&\cdots&0&0\\
                     -1&0&\cdots&0&0\\
                     \vdots&\vdots&\ddots&\vdots&\vdots\\
                     0&0&\cdots&0&1\\
                     0&0&\cdots&-1&0
\end{matrix}\right).
\een
The images
of the elements $\si_{ij}$
under the isomorphism \eqref{speci}
are respectively the elements $F_{ij}$ of $\spa_{2n}$.

In the next section we prove that if $q$ is specialized to a nonzero complex
number such that $q^2\ne 1$, then the
mapping \eqref{atemb} defines an embedding
of the respective algebras over $\CC$.

\section{Poincar\'e--Birkhoff--Witt theorem}\label{sec:pbw}
\setcounter{equation}{0}

Define the extended twisted
quantized enveloping algebra $\check\U'_q(\spa_{2n})$
as follows.
This is an associative algebra
generated by elements $s_{ij}$,
$i,j\in\{1,\dots,2n\}$ and $s_{i,i+1}^{-1}$, $i=1,3,\dots,2n-1$,
where
$s_{ij}=0$
for $j=i+1$ with even $i$, and
for $j\geqslant i+2$ and all $i$.
The defining relations are given by
\eqref{invrel} and \eqref{drabss}.
As with the algebra $\U'_q(\spa_{2n})$,
we combine the $s_{ij}$ into
a matrix $S$ which a block-triangular form with
$n$ diagonal $2\times 2$-blocks.

Consider the algebra $\hat\U_q(\gl_{2n})$ introduced in the previous
section, and denote by $\check\U_q(\gl_{2n})$ its quotient
by the ideal generated by the central elements $t_{ii}\tss\bar t_{ii}-1$
for all even $i$. We keep the same notation
for the images of the generators of $\hat\U_q(\gl_{2n})$
in $\check\U_q(\gl_{2n})$.
The mapping given by
\beql{atembext}
\psi:S\mapsto T\ts G\ts \overline T^{\ts t}
\eeq
defines a homomorphism
$\check\U'_q(\spa_{2n})\to\check\U_q(\gl_{2n})$.
This is verified by the same calculation as with the homomorphism
\eqref{atemb}; see \cite{n:ms} and \cite{mrs:cs}.
In particular,
\ben
s_{i,i+1}\mapsto q\ts t_{ii}\ts\bar t_{i+1,i+1}\Fand
s^{-1}_{i,i+1}\mapsto q^{-1}\ts t^{-1}_{ii}\ts t^{}_{i+1,i+1}
\een
for $i=1,3,\dots,2n-1$.

\bth\label{thm:embext}
The mapping \eqref{atembext} defines an embedding
$\check\U'_q(\spa_{2n})\hra\check\U_q(\gl_{2n})$.
\eth

\bpf
We only need to show that the kernel of the homomorphism \eqref{atembext}
is zero. We shall use a weak form of the Poincar\'e--Birkhoff--Witt theorem
for the algebra $\check\U'_q(\spa_{2n})$ provided by the
following lemma.

\ble\label{lem:wpbw}
The algebra $\check\U'_q(\spa_{2n})$ is spanned by
the monomials
\begin{multline}\label{ordmonext}
{s}_{2n,1}^{\ts k_{2n,1}}\dots \ts {s}_{2n,2n}^{\ts k_{2n,2n}}\ts
{s}_{2n-2,1}^{\ts k_{2n-2,1}}\dots \ts {s}_{2n-2,2n-2}^{\ts k_{2n-2,2n-2}}\ts
\dots\ts
{s}_{21}^{\ts k_{21}}\ts {s}_{22}^{\ts k_{22}}\\
{}\times
{s}_{2n-1,1}^{\ts k_{2n-1,1}}\dots \ts {s}_{2n-1,2n}^{\ts k_{2n-1,2n}}\ts
{s}_{2n-3,1}^{\ts k_{2n-3,1}}\dots \ts {s}_{2n-3,2n-2}^{\ts k_{2n-3,2n-2}}\ts
{s}_{11}^{\ts k_{11}}\ts {s}_{12}^{\ts k_{12}},
\end{multline}
where $k_{12}$, $k_{34},\dots,k_{2n-1,2n}$ run over all integers  while
the remaining $k_{ij}$ run over non-negative integers.
\ele

\bpf
We follow the argument of \cite{mrs:cs}, where a similar statement
was proved (see Lemma~2.1 there).
We shall be proving that any monomial
in the generators can be written as a linear combination of
monomials of the form \eqref{ordmonext}.
The defining relations \eqref{drabss} imply that
\beql{soopo}
q^{-\de_{ik}+\de_{i+1,k}}\ts
s^{}_{i,i+1}\ts s^{}_{kl}=q^{\de_{il}-\de_{i+1,l}}\ts
s^{}_{kl}\ts s^{}_{i,i+1},
\eeq
for any $i=1,3,\dots,2n-1$. Therefore, it suffices to consider
monomials where the generators $s_{12},\dots,s_{2n-1,2n}$
occur with non-negative powers.
For any monomial
\beql{anymon}
s_{i_1a_1}\cdots\ts s_{i_pa_p}
\eeq
we introduce its {\it weight\/} $w$ by $w=i_1+\cdots+i_p$.
We shall use induction on $w$.
The defining relations \eqref{drabss} for $\check\U'_q(\spa_{2n})$
imply that, modulo products of weight less than $i+j$, we have
\beql{igj}
q^{\delta_{aj}+\delta_{ij}}\ts s_{ia}\ts s_{jb}\equiv
q^{\delta_{ab}+\delta_{ib}}\ts s_{jb}\ts s_{ia}
+(q-\qin)\ts q^{\delta_{ai}}\ts(\de_{b<a} -\de_{i<j})
\ts s_{ja}\ts s_{ib}.
\eeq
Let $\pi$ denote the permutation of the indices $1,2,\dots,2n$ such that
$\pi^{-1}$ is given by $(2n,2n-2,\dots,2,2n-1,2n-3,\dots,1)$. Relation \eqref{igj}
allows us to represent \eqref{anymon}, modulo
monomials of weight less than $w$,
as a linear combination of monomials
$s_{j_1b_1}\cdots\ts s_{j_pb_p}$ of weight $w$ such that
$\pi(j_1)\leqslant\cdots \leqslant \pi(j_p)$. Consider a sub-monomial
$s_{i\tss c_1}\cdots\ts s_{i\tss c_r}$
containing generators with the same first index.
By \eqref{igj} we have
\beql{igi}
s_{ia}\ts s_{ib}\equiv
q^{\delta_{ib}-\delta_{ia}+1}\ts s_{ib}\ts s_{ia}
\end{equation}
for $a>b$. Using this relation repeatedly we
bring the sub-monomial to the required form.
\epf

\ble\label{lem:linind}
The monomials
\begin{multline}\label{ordmonextdo}
{s}_{2n,1}^{\ts k_{2n,1}}\dots \ts {s}_{2n,2n-2}^{\ts k_{2n,2n-2}}
\ts {s}_{2n,2n}^{\ts k_{2n,2n}}\dots\ts
{s}_{41}^{\ts k_{41}}\ts {s}_{42}^{\ts k_{42}}\ts {s}_{44}^{\ts k_{44}}
\ts {s}_{22}^{\ts k_{22}}\\
{}\times
{s}_{2n-1,1}^{\ts k_{2n-1,1}}\dots \ts {s}_{2n-1,2n-2}^{\ts k_{2n-1,2n-2}}\ts
\dots\ts {s}_{31}^{\ts k_{31}}\ts {s}_{32}^{\ts k_{32}}\\
{}\times
{s}_{2n,2n-1}^{\ts k_{2n,2n-1}}\ts {s}_{2n-1,2n}^{\ts k_{2n-1,2n}}
\dots\ts
{s}_{21}^{\ts k_{21}}\ts {s}_{12}^{\ts k_{12}}\ts
{s}_{2n-1,2n-1}^{\ts k_{2n-1,2n-1}}\dots\ts
{s}_{11}^{\ts k_{11}},
\end{multline}
where $k_{12}$, $k_{34},\dots,k_{2n-1,2n}$ run over all integers  while
the remaining $k_{ij}$ run over non-negative integers,
are linearly independent in $\check\U'_q(\spa_{2n})$.
\ele

\bpf
Let $\mu_i$ and $\bar\mu_i$ with $i=1,\dots,2n$ be arbitrary
nonzero complex numbers
such that $\mu_i=\bar\mu_i=1$ for all even $i$.
Consider the corresponding Verma module
$M(\mu,\bar\mu)$ over the algebra $\check\U_q(\gl_{2n})$
which is defined as the quotient of $\check\U_q(\gl_{2n})$
by the left ideal generated by the elements $\bar t_{ij}$ with $i<j$
and $t_{ii}-\mu_i$, $\bar t_{ii}-\bar\mu_i$ with $i=1,\dots,2n$.
Corollary~\ref{cor:extpbw} implies that the elements
\ben
t_{i_1a_1}^{k_{i_1a_1}}\dots\ts
t_{i_ma_m}^{k_{i_ma_m}}\ts \xi
\een
form a basis of the Verma module $M(\mu,\bar\mu)$,
where $\xi$ denotes its highest vector
and the generators $t_{i_ra_r}$ with $i_r>a_r$ are written in
accordance with a certain linear ordering determined by the
permutation $\pi$ defined in the proof of Lemma~\ref{lem:wpbw}.

Suppose now that a nontrivial linear combination
of monomials \eqref{ordmonextdo} is zero.
For any odd $i$ the image of the power $s_{ii}^k$ under the homomorphism
\eqref{atembext} is given by
\beql{imsii}
\psi:s_{ii}^k\mapsto q^{\frac{3k-k^2}{2}}
\ts t_{ii}^{\ts k}\ts \bar t_{i,i+1}^{\ts k}.
\eeq
Amongst the monomials which occur in the linear
combination, pick up
a monomial such that each power
$k_{ii}$ takes a minimal possible value $\kappa_i$ for each odd $i$.
Now take the image of the linear combination
under \eqref{atembext} and apply this image
to the vector
\ben
t_{21}^{\kappa_1}\dots t_{2n,2n-1}^{\kappa_{2n-1}}\ts\xi\in M(\mu,\bar\mu).
\een
By the choice the parameters $\kappa_i$, the nonzero
contribution to the resulting expression can only come from
the monomials \eqref{ordmonextdo} with $k_{ii}=\kappa_i$ for all odd $i$.
However, by \eqref{imsii}
\ben
\psi(s_{ii}^{\kappa_i})\ts t_{i+1,i}^{\kappa_i}\ts\xi=c_i\ts \xi,
\een
where $c_i$ is a constant depending on the
$\mu_j$ and $\bar\mu_j$ which can easily be calculated.
Choosing the parameters $\mu_j$ and $\bar\mu_j$ in such
a way that $c_i\ne 0$ for each odd $i$,
we conclude
that the image under $\psi$ of a certain nontrivial linear combination of monomials
of the form
\begin{multline}\label{ordmonextdowi}
{s}_{2n,1}^{\ts k_{2n,1}}\dots \ts {s}_{2n,2n-2}^{\ts k_{2n,2n-2}}
\ts {s}_{2n,2n}^{\ts k_{2n,2n}}\dots\ts
{s}_{41}^{\ts k_{41}}\ts {s}_{42}^{\ts k_{42}}\ts {s}_{44}^{\ts k_{44}}
\ts {s}_{22}^{\ts k_{22}}\\
{}\times
{s}_{2n-1,1}^{\ts k_{2n-1,1}}\dots \ts {s}_{2n-1,2n-2}^{\ts k_{2n-1,2n-2}}\ts
\dots\ts {s}_{31}^{\ts k_{31}}\ts {s}_{32}^{\ts k_{32}}\\
{}\times
{s}_{2n,2n-1}^{\ts k_{2n,2n-1}}\ts {s}_{2n-1,2n}^{\ts k_{2n-1,2n}}
\dots\ts
{s}_{21}^{\ts k_{21}}\ts {s}_{12}^{\ts k_{12}}
\end{multline}
acts as zero when applied to the highest vector $\xi$ of the Verma module
$M(\mu,\bar\mu)$. As the image of $s_{ij}$ under the homomorphism $\psi$
is given by \eqref{sijtij}, we find that
\ben
\psi(s_{i,i+1})\ts\xi=q\ts \mu_i\ts\xi,\qquad
\psi(s_{i+1,i})\ts\xi=-\bar\mu_{i}\ts\xi
\een
for any odd $i$. Moreover, using \eqref{sijtij} again, we come to the formulas
\ben
\psi(s_{ij})\ts\xi=\begin{cases}
                   q\ts t_{i,j-1}\ts\xi,
                   \qquad&\text{if $j$ is even, $i\geqslant j$},\\
                   -\bar\mu_j\ts t_{i,j+1}\ts\xi,
                   \qquad&\text{if $j$ is odd, $i\geqslant j+2$}.
                   \end{cases}
\een
Varying the parameters $\mu_i$ and $\bar\mu_i$,
we conclude that the elements
\begin{multline}\label{tordmonewiip}
{t}_{2n,2}^{\ts k_{2n,1}}\dots \ts {t}_{2n,2n-3}^{\ts k_{2n,2n-2}}
\ts {t}_{2n,2n-1}^{\ts k_{2n,2n}}\dots\ts
{t}_{42}^{\ts k_{41}}\ts {t}_{41}^{\ts k_{42}}\ts {t}_{43}^{\ts k_{44}}
\ts {t}_{21}^{\ts k_{22}}\\
{}\times
{t}_{2n-1,2}^{\ts k_{2n-1,1}}\dots \ts {t}_{2n-1,2n-3}^{\ts k_{2n-1,2n-2}}\ts
\dots\ts {t}_{32}^{\ts k_{31}}\ts {t}_{31}^{\ts k_{32}}
\end{multline}
of the algebra $\check\U_q(\gl_{2n})$ are linearly dependent.
This contradicts Corollary~\ref{cor:extpbw}, completing the proof
of the lemma.
\epf

Now, let us denote by $\check\U^{\prime\tss +}$ the subalgebra
of $\check\U'_q(\spa_{2n})$ generated by the elements
$s_{i,i+1}$ for all odd $i$ and $s_{ij}$ for all $i\geqslant j$.
For any non-negative integer $m$
consider the subspace $\check\U^{\prime\tss +}_m$ of
$\check\U^{\prime\tss +}$
of elements of degree at most $m$ in the generators.
Lemma~\ref{lem:wpbw} implies that
$\check\U^{\prime\tss +}_m$ is spanned by the monomials
\eqref{ordmonext} of total degree ${}\leqslant m$
with non-negative powers of the generators.
On the other hand, by Lemma~\ref{lem:linind},
the monomials \eqref{ordmonextdo} of total degree ${}\leqslant m$
with non-negative powers of the generators
are linearly independent.
Since the numbers of both types of monomials
of total degree ${}\leqslant m$ coincide,
we conclude that each of these families of monomials forms
a basis of $\check\U^{\prime\tss +}_m$. This implies that
each family of monomials \eqref{ordmonext}
and \eqref{ordmonextdo} forms a basis of $\check\U'_q(\spa_{2n})$.
In particular, any element $u$ of $\check\U'_q(\spa_{2n})$
can be written as linear combination of monomials
\eqref{ordmonextdo}. However, the proof of Lemma~\ref{lem:linind}
shows that $\psi(u)=0$ implies $u=0$, thus proving
that the kernel of $\psi$ is zero.
\epf

The following
version of the
Poincar\'e--Birkhoff--Witt theorem for the algebra
$\check\U'_q(\spa_{2n})$ was already noted
in the proof of Theorem~\ref{thm:embext}.

\bco\label{cor:pbw}
Each family of monomials \eqref{ordmonext}
and \eqref{ordmonextdo} constitutes a basis of
the algebra $\check\U'_q(\spa_{2n})$.
\qed
\eco

For any $i=1,3,\dots,2n-1$ set
\ben
\vartheta_i=
s^{}_{i+1,i+1}\ts s^{}_{ii}-q^2\ts s^{}_{i+1,i}s^{}_{i,i+1}.
\een
As was observed in \cite[Section~2.2]{mrs:cs}, the elements
$\vartheta_{i}$ belong to the center of $\check\U'_q(\spa_{2n})$.

\bpr\label{prop:basthe}
The elements
\begin{multline}\label{ordmonextth}
{s}_{2n,1}^{\ts k_{2n,1}}\dots \ts {s}_{2n,2n-2}^{\ts k_{2n,2n-2}}\ts
{s}_{2n,2n}^{\ts k_{2n,2n}}\ts \vartheta_{2n-1}^{\ts k_{2n-1}}\ts
{s}_{2n-1,2n}^{\ts k_{2n-1,2n}}\ts{s}_{2n-1,2n-1}^{\ts k_{2n-1,2n-1}}\\
{}\times\dots\times {s}_{41}^{\ts k_{41}}\ts {s}_{42}^{\ts k_{42}}
\ts {s}_{44}^{\ts k_{44}}\ts \vartheta_{3}^{\ts k_{3}}
\ts {s}_{34}^{\ts k_{34}}\ts {s}_{33}^{\ts k_{33}}
\ts {s}_{22}^{\ts k_{22}}\ts\vartheta_{1}^{\ts k_{1}}
\ts {s}_{12}^{\ts k_{12}}\ts{s}_{11}^{\ts k_{11}}\\
{}\times {s}_{31}^{\ts k_{31}}\ts {s}_{32}^{\ts k_{32}}\dots\ts
{s}_{2n-1,1}^{\ts k_{2n-1,1}}\dots \ts
{s}_{2n-1,2n-2}^{\ts k_{2n-1,2n-2}},
\end{multline}
where the $k_{i,i+1}$ for odd $i$ run over all integers while $k_i$
and the remaining $k_{rj}$ run over non-negative
integers, form a basis of the algebra $\check\U'_q(\spa_{2n})$.
\epr

\bpf
First, we prove that a basis of $\check\U'_q(\spa_{2n})$ is comprised by
the monomials
\begin{multline}\label{ordmnb}
{s}_{2n,1}^{\ts k_{2n,1}}\dots \ts {s}_{2n,2n-2}^{\ts k_{2n,2n-2}}\ts
{s}_{2n,2n}^{\ts k_{2n,2n}}\ts {s}_{2n,2n-1}^{\ts k_{2n,2n-1}}\ts
{s}_{2n-1,2n}^{\ts k_{2n-1,2n}}\ts{s}_{2n-1,2n-1}^{\ts k_{2n-1,2n-1}}\\
{}\times\dots\times {s}_{41}^{\ts k_{41}}\ts {s}_{42}^{\ts k_{42}}
\ts {s}_{44}^{\ts k_{44}}\ts {s}_{43}^{\ts k_{43}}
\ts {s}_{34}^{\ts k_{34}}\ts {s}_{33}^{\ts k_{33}}
\ts {s}_{22}^{\ts k_{22}}\ts{s}_{21}^{\ts k_{21}}
\ts {s}_{12}^{\ts k_{12}}\ts{s}_{11}^{\ts k_{11}}\\
{}\times {s}_{31}^{\ts k_{31}}\ts {s}_{32}^{\ts k_{32}}\dots\ts
{s}_{2n-1,1}^{\ts k_{2n-1,1}}\dots \ts
{s}_{2n-1,2n-2}^{\ts k_{2n-1,2n-2}},
\end{multline}
where the $k_{i,i+1}$ for odd $i$ are integers
and the remaining $k_{rj}$ are non-negative
integers. Indeed, the argument used in the proof of Lemma~\ref{lem:wpbw}
together with \eqref{soopo}
implies that monomials \eqref{ordmnb} span the algebra
$\check\U'_q(\spa_{2n})$; the only additional observation required
is that
\beql{siiodd}
s_{ii}\ts s_{jb}=s_{jb}\ts s_{ii}
\eeq
for all odd $i$ such that $i\geqslant j\geqslant b$; see \eqref{drabss}.
Then, as in the proof of Theorem~\ref{thm:embext}, we conclude
that monomials \eqref{ordmnb} form a basis of $\check\U'_q(\spa_{2n})$.

Now we show that the elements \eqref{ordmonextth}
span the algebra $\check\U'_q(\spa_{2n})$.
Since
the defining relations between the generators
$s_{ii}$, $s_{i,i+1}$, $s_{i+1,i}$ and $s_{i+1,i+1}$
do not involve any other generators,
it is sufficient to consider the particular case $n=1$.
We have
\ben
{s}_{22}^{\ts k}\ts s^{\ts l}_{21}
\ts {s}_{12}^{\ts m}\ts{s}_{11}^{\ts r}
=-q^{-2}\ts {s}_{22}^{\ts k}\ts \vartheta^{}_{1}
\tss s^{\ts l-1}_{21}
\tss {s}_{12}^{\ts m-1}\tss{s}_{11}^{\ts r}+\sum_{a=1}^{l}
c_a\ts {s}_{22}^{\ts k+1}
\tss s^{\ts a-1}_{21}
\tss {s}_{12}^{\ts m+l-a-1}\tss{s}_{11}^{\ts r+1}
\een
for some complex coefficients $c_a$. Arguing by the
induction on $l$, we find that
\begin{multline}\label{relsmde}
{s}_{22}^{\ts k}\ts s^{\ts l}_{21}
\ts {s}_{12}^{\ts m}\ts{s}_{11}^{\ts r}=(-1)^l\ts q^{-2l}
\ts {s}_{22}^{\ts k}\ts \vartheta^{\ts l}_{1}
\tss {s}_{12}^{\ts m-l}\tss{s}_{11}^{\ts r}\\
+\quad\text{a linear combination of\quad
${s}_{22}^{\ts k'}\ts \vartheta^{\ts l'}_{1}
\tss {s}_{12}^{\ts m'}\tss{s}_{11}^{\ts r'}$\quad with\quad $l'<l$}.
\end{multline}
Thus, the elements \eqref{ordmonextth}
span $\check\U'_q(\spa_{2n})$.
The relations \eqref{relsmde} can obviously be inverted
to get similar expressions of the elements
${s}_{22}^{\ts k}\ts \vartheta^{\ts l}_{1}
\tss {s}_{12}^{\ts m}\tss{s}_{11}^{\ts r}$
in terms of monomials ${s}_{22}^{\ts k}\ts s^{\ts l}_{21}
\ts {s}_{12}^{\ts m}\ts{s}_{11}^{\ts r}$.
This implies that the elements \eqref{ordmonextth}
are linearly independent.
\epf

By the definition of the algebra $\U'_q(\spa_{2n})$,
we have a surjective homomorphism
$\check\U'_q(\spa_{2n})\to\U'_q(\spa_{2n})$ which takes the
generators $s_{ij}$ to the elements of $\U'_q(\spa_{2n})$
with the same name. In other words, we have an isomorphism
\ben
\U'_q(\spa_{2n})\cong \check\U'_q(\spa_{2n})/I,
\een
where $I$ is the ideal of $\check\U'_q(\spa_{2n})$ generated
by the central elements $\vartheta_i-q^3$ for all odd
$i=1,3,\dots,2n-1$.
Hence, the following Poincar\'e--Birkhoff--Witt theorem for the algebra
$\U'_q(\spa_{2n})$ follows from Proposition~\ref{prop:basthe}
and the relations \eqref{soopo} and \eqref{siiodd}.

\bth\label{thm:pbw}
The elements
\begin{multline}\label{ordpbw}
{s}_{2n,1}^{\ts k_{2n,1}}\dots \ts {s}_{2n,2n-2}^{\ts k_{2n,2n-2}}\ts
{s}_{2n,2n}^{\ts k_{2n,2n}}\dots\ts
{s}_{41}^{\ts k_{41}}\ts {s}_{42}^{\ts k_{42}}
\ts {s}_{44}^{\ts k_{44}}\ts {s}_{22}^{\ts k_{22}}\\
{}\times
{s}_{11}^{\ts k_{11}}\ts {s}_{12}^{\ts k_{12}}\ts
{s}_{31}^{\ts k_{31}}\ts {s}_{32}^{\ts k_{32}}\ts {s}_{33}^{\ts k_{33}}
\ts {s}_{34}^{\ts k_{34}}\dots\ts
{s}_{2n-1,1}^{\ts k_{2n-1,1}}\dots \ts
{s}_{2n-1,2n}^{\ts k_{2n-1,2n}},
\end{multline}
where the $k_{i,i+1}$ for odd $i$ run over all integers
and the remaining $k_{rj}$ run over non-negative
integers, form a basis of the algebra $\U'_q(\spa_{2n})$.
\qed
\eth

\bco\label{cor:embquo}
The mapping \eqref{atembext} defines an
embedding $\U'_q(\spa_{2n})\hra\U_q(\gl_{2n})$.
\eco

\bpf
By Theorem~\ref{thm:embext}, the algebra $\check\U'_q(\spa_{2n})$
can be identified with a subalgebra of $\check\U_q(\gl_{2n})$.
Then, for each $i=1,3,\dots,2n-1$ the relation
$\vartheta_i=q^3$ is equivalent to $t_{ii}\tss\bar t_{ii}=1$.
However, the quotient of $\check\U'_q(\spa_{2n})$ by
the relations $t_{ii}\tss\bar t_{ii}=1$ for all odd $i$
is isomorphic to $\U_q(\gl_{2n})$. Hence, the claim follows
from Theorem~\ref{thm:pbw}.
\epf

Using Corollary~\ref{cor:embquo}, we shall regard
$\U'_q(\spa_{2n})$ as a subalgebra of $\U_q(\gl_{2n})$.

\bre\label{rem:osyquo}
An analogue of the Poincar\'e--Birkhoff--Witt theorem for the algebra
$\U'_q(\oa_{N})$ was proved in \cite{ik:nd}
with the use of the Diamond Lemma. If
$\U'_q(\oa_{N})$ is regarded as an algebra over $\CC(q)$,
the theorem can also be proved by a specialization argument;
see \cite[Corollary~2.3]{mrs:cs}.
If $q$ is a nonzero complex number such that $q^2\ne 1$, a proof
can be given in a way similar
to the above arguments with some simplifications. Indeed,
recall that $\U'_q(\oa_{N})$ is generated by elements
$s_{ij}$ with the defining relations
\begin{align}
s_{ij}&=0, \qquad 1 \leq i<j\leqslant N,
\non\\
s_{ii}&=1,\qquad 1\leqslant i\leqslant N,
\non\\
R\ts S_1& R^{\ts\prime} S_2=S_2R^{\ts\prime} S_1R,
\non
\end{align}
using the same notation as for the symplectic case.
The monomials
\beql{monorth}
{s}_{21}^{\ts k_{21}}\ts    {s}_{31}^{\ts k_{31}}\ts
{s}_{32}^{\ts k_{32}}\ts  \cdots \ts
{s}_{N1}^{\ts k_{N1}}\ts {s}_{N2}^{\ts k_{N2}} \ts
\cdots\ts  {s}_{N,N-1}^{\ts k_{N,N-1}}
\eeq
span the algebra $\U'_q(\oa_{N})$; see \cite[Lemma~2.1]{mrs:cs}.
In order to prove the linear independence of these
monomials, consider their images under the
homomorphism $\phi:\U'_q(\oa_{N})\to\U_q(\gl_N)$
defined by $S\mapsto T\ts \overline T^{\ts t}$.
Consider the Verma module $M(\mu)$ over $\U_q(\gl_N)$
with $\mu=(1,\dots,1)$ and the highest vector $\xi$ so that
\ben
\bar t_{ij}\ts\xi=0,\qquad i<j\Fand t_{ii}\ts\xi=1,\qquad i=1,\dots,N.
\een
Applying the image of the monomial \eqref{monorth} to
the highest vector $\xi$ we get
\ben
{t}_{21}^{\ts k_{21}}\ts    {t}_{31}^{\ts k_{31}}\ts
{t}_{32}^{\ts k_{32}}\ts  \cdots \ts
{t}_{N1}^{\ts k_{N1}}\ts {t}_{N2}^{\ts k_{N2}} \ts
\cdots\ts  {t}_{N,N-1}^{\ts k_{N,N-1}}\ts \xi.
\een
The proof is completed by using Proposition~\ref{prop:arbord}.
This argument also shows that the homomorphism $\phi$
is an embedding.
\qed
\ere

\section{Highest weight representations}\label{sec:hw}
\setcounter{equation}{0}

Here we introduce highest weight representations of $\U'_q(\spa_{2n})$ and
prove that every finite-dimensional irreducible representation
of this algebra is highest weight. The arguments
are quite standard; cf. \cite{h:il}, \cite[Chapter~10]{cp:gq}.
From now on, $q$ will denote a complex number which is nonzero
and not a root of unity.

We shall be using the following notation.
For any two $n$-tuples $\al=(\al_1,\dots,\al_n)$ and
$\be=(\be_1,\dots,\be_n)$
we shall denote by $\al\cdot\be$ the $n$-tuple
$(\al_1\be_1,\dots,\al_n\be_n)$.
Also, $q^{\tss\al}$ will denote the $n$-tuple
$(q^{\tss\al_1},\dots,q^{\tss\al_n})$ in the case where $\al$
is an $n$-tuple of integers.

A representation $V$ of $\U'_q(\spa_{2n})$ will be called
a {\it highest weight representation\/} if $V$ is generated
by a nonzero vector $v$ such that
\begin{alignat}{2}
s_{ij}\ts v&=0 \qquad &&\text{for} \quad i=1,3,\dots,2n-1,\quad
j=1,2,\dots,i, \qquad \text{and}
\non\\
s_{i,i+1}\ts v&=\la_i\ts v\qquad &&\text{for} \quad i=1,3,\dots,2n-1,
\non
\end{alignat}
for some complex numbers $\la_i$. These numbers have to be nonzero due to
the relation \eqref{invrel}.
The $n$-tuple $\la=(\la_1,\la_3,\dots,\la_{2n-1})$ will be called
the {\it highest weight\/} of $V$.
Observe that due to the relations \eqref{drabss},
the elements $s^{}_{i,i+1}$ with odd $i$ pairwise commute.
The commutative subalgebra of $\U'_q(\spa_{2n})$ generated
by these elements will play the role of a Cartan subalgebra.

Consider the root system $\De$ of type $C_n$ which is the subset
of vectors in $\RR^n$ of the form
\ben
{}\pm 2\ts\ve_i\quad\text{with}\quad 1\leqslant i\leqslant n\Fand
{}\pm\ve_i\pm\ve_j\quad\text{with}\quad 1\leqslant i<j\leqslant n,
\een
where $\ve_i$ denotes the $n$-tuple which has $1$ on the $i$-th position
and zeros elsewhere. Partition this set into positive
and negative roots $\De=\De^+\cup(-\De^+)$, where the set
of positive roots $\De^+$ consists of the vectors
\ben
2\ts\ve_i\quad\text{with}\quad 1\leqslant i\leqslant n\Fand
\ve_i+\ve_j,\quad {}-\ve_i+\ve_j\quad\text{with}
\quad 1\leqslant i<j\leqslant n.
\een

For any $n$-tuple of nonzero complex numbers
$\mu=(\mu_1,\mu_3,\dots,\mu_{2n-1})$ define the corresponding
{\it weight subspace\/} of $V$ by
\ben
V_{\mu}=\{w\in V\ |\ s_{i,i+1}\ts w=\mu_i\ts w\qquad
\text{for} \quad i=1,3,\dots,2n-1\}.
\een
Any nonzero vector $w\in V_{\mu}$
is called a {\it weight vector\/} of weight
$\mu$.

Given an $n$-tuple of nonzero complex numbers
$\la=(\la_1,\la_3,\dots,\la_{2n-1})$,
the corresponding {\it Verma module\/} $M(\la)$ over $\U'_q(\spa_{2n})$
is defined as the quotient of $\U'_q(\spa_{2n})$ by the left ideal
generated by the elements
\beql{uptri}
s_{ij}\quad\text{with} \quad i=1,3,\dots,2n-1,\quad
j=1,2,\dots,i,
\eeq
and
\ben
s_{i,i+1}-\la_i\quad\text{with} \quad i=1,3,\dots,2n-1.
\een
The Verma module is obviously a highest weight representation.
The image $\xi$ of the element $1\in \U'_q(\spa_{2n})$
in $M(\la)$ is the highest vector and $\la$ is the highest weight.
By the Poincar\'e--Birkhoff--Witt theorem for the algebra
$\U'_q(\spa_{2n})$ (see Theorem~\ref{thm:pbw}), a basis of $M(\la)$ is comprised by
the elements
\ben
\prod_{i=2,4,\dots,2n}^{\leftarrow}
s_{i1}^{k_{i1}}\ts s_{i2}^{k_{i2}}\ts
\dots\ts s_{i,i-2}^{k_{i,i-2}}\ts s_{i,i}^{k_{i,i}}\ts \xi,
\een
where the $k_{kj}$ run over non-negative integers.
Due to \eqref{soopo},
we have the weight space decomposition
\ben
M(\la)=\underset{\mu}{\oplus}\ts M(\la)_{\mu}.
\een
The weight subspace $M(\la)_{\mu}$ is nonzero if and only if
$\mu$ has the form $\mu=q^{-\om}\cdot\la$, where
$\om$ is a linear combination of elements of $\De^+$ with non-negative
integer coefficients. The dimension of $M(\la)_{\mu}$
is given by the same formula
as in the classical case; see e.g. \cite{h:il}.
In particular, the weight space $M(\la)_{\la}$ is one-dimensional
and spanned by $\xi$.

Every highest weight module $V$ with the highest weight $\la$
is a homomorphic image of $M(\la)$. So, $V$ is
the direct sum of its weight subspaces,
$
V=\oplus\ts V_{\mu}.
$
By a standard argument, $M(\la)$ contains a unique maximal submodule
which does not contain the vector $\xi$. The quotient of $M(\la)$
by this submodule is, up to an isomorphism,
the unique irreducible highest weight module with the highest weight $\la$.
We shall denote this quotient by $L(\la)$.

\bpr\label{prop:hwpro}
Every finite-dimensional irreducible representation $V$ of
$\ts\U'_q(\spa_{2n})$ is isomorphic to $L(\la)$ for a certain
highest weight $\la$.
\epr

\bpf
This is verified by a standard argument.
We need to show
that $V$ contains a weight vector annihilated by all
operators \eqref{uptri}.
Since the operators $s_{i,i+1}$ with $i=1,3,\dots,2n-1$
on $V$
pairwise commute, the space $V$ must contain
a weight vector $w$ of a certain weight $\mu$.
If $w$ is not annihilated by the operators \eqref{uptri} then
applying these operators to $w$
we can get other weight vectors with weights of the form
$q^{\om}\cdot \mu$, where $\om$ is
a linear combination of elements of $\De^+$ with non-negative
integer coefficients. As $\dim V<\infty$, the proof is completed
by the classical argument; see e.g. \cite{h:il}.
\epf

Due to Proposition~\ref{prop:hwpro}, in order to describe
the finite-dimensional irreducible representations of $\U'_q(\spa_{2n})$
it suffices to find the necessary and sufficient conditions on $\la$
for the representation $L(\la)$ to be finite-dimensional.
As in the classical theory, the case $n=1$ plays
an important role.

\section{Representations of $\U'_q(\spa_{2})$}\label{sec:sptwo}
\setcounter{equation}{0}

Using \eqref{qdetrel}, we can regard $\U'_q(\spa_{2})$
as the algebra with generators $s^{}_{11},s^{}_{22},s^{}_{12},s^{-1}_{12}$.
The defining relations take the form of \eqref{invrel} with $i=1$
together with
\beql{drspt}
s^{}_{11}\tss s^{}_{22}=q^{-2} s^{}_{22}\tss s^{}_{11}-(q-\qin)(s^2_{12}-q^2)
\eeq
and
\beql{sonetwo}
s^{}_{12}\tss s^{}_{11}=q^{2}\tss s^{}_{11}\tss s^{}_{12},\qquad
s^{}_{12}\tss s^{}_{22}=q^{-2}\tss s^{}_{22}\tss s^{}_{12}.
\eeq
For a nonzero complex number $\la$, the corresponding Verma module
$M(\la)$ has the basis
\ben
s_{22}^k\ts \xi,\qquad k\geqslant 0.
\een
Using \eqref{drspt} and \eqref{sonetwo} we obtain
\ben
s^{}_{11}\tss s_{22}^k\ts \xi=q^3\tss (1-\la^2\tss
q^{-2k})(1-q^{-2k})\tss s_{22}^{k-1}\ts \xi.
\een
Hence, the module $M(\la)$ is reducible if and only if $\la=\si\tss q^m$
for a positive integer $m$ and $\si\in\{-1,1\}$.
Thus, we have the following.

\bpr\label{prop:sltwoir}
The irreducible highest weight module $L(\la)$
over $\U'_q(\spa_{2})$ is finite-dimensional
if and only if $\la=\si\tss q^m$
for a positive integer $m$ and $\si\in\{-1,1\}$.
\qed
\epr

In this case $L(\la)$ has a basis
$\{\tss v^{}_k\tss\}$, $k=0,1,\dots,m-1$, with
the action of $\U'_q(\spa_{2})$ given by
\ben
s^{}_{12}\tss v^{}_k=\si\tss q^{m-2k}\tss v^{}_k,
\een
\ben
s^{}_{22}\tss v^{}_k=v^{}_{k+1},
\een
\ben
s^{}_{11}\tss v^{}_k=q^3\tss (1-q^{2m-2k})(1-q^{-2k})\tss v^{}_{k-1},
\een
where $v^{}_{-1}=v^{}_{m}=0$.

Note also that
all finite-dimensional irreducible $\U'_q(\spa_{2})$-modules
can be obtained by restriction from $\U_q(\gl_2)$-modules.
Indeed, consider the irreducible highest weight module
over $\U_q(\gl_2)$ generated by a vector $w$ satisfying
\ben
\bar t_{12}\ts w=0,\qquad t_{11}\ts w=\vs_1\tss q^{\tss\mu_1}\ts w,
\qquad t_{22}\ts w=\vs_2\tss q^{\tss\mu_2}\ts w,
\een
where $\vs_1,\vs_2\in\{-1,1\}$ and $\mu_1,\mu_2\in\ZZ$.
If $\mu_1-\mu_2\geqslant 0$ then this module
has dimension $m=\mu_1-\mu_2+1$ and its restriction to the subalgebra
$\U'_q(\spa_{2})$ is isomorphic to $L(\la)$ with
$\la=\si\tss q^m$, where $\si=\vs_1/\vs_2$.
This is easily derived by using \eqref{sijtij}.

\section{Classification theorem}\label{sec:cth}
\setcounter{equation}{0}

Consider an arbitrary irreducible highest weight
representation $L(\la)$ of $\U'_q(\spa_{2n})$.
Our first aim is to find necessary conditions
on $\la=(\la_1,\la_3,\dots,\la_{2n-1})$ for $L(\la)$ to be finite-dimensional.
So, suppose that $\dim L(\la)<\infty$.

For any index $k=1,2,\dots,n$ the subalgebra of $\U'_q(\spa_{2n})$
generated by the elements $s^{}_{2k-1,2k-1}$, $s^{}_{2k,\ts 2k}$,
$s^{}_{2k-1,2k}$ and
$s^{-1}_{2k-1,2k}$ is isomorphic to $\U'_q(\spa_{2})$.
The cyclic span of the highest vector $\xi$ of $L(\la)$
with respect to this subalgebra is a highest weight module
with the highest weight $\la_{2k-1}$. The irreducible quotient
of this module is finite-dimensional and so, by Proposition~\ref{prop:sltwoir},
we must have $\la_{2k-1}=\si_k\ts q^{m_k}$ for some
positive integer $m_k$ and $\si_k\in\{-1,1\}$. Thus, the highest
weight $\la$ of $L(\la)$ must have the form
\ben
\la=(\si_1\ts q^{m_1},\dots,\si_n\ts q^{m_n}).
\een

Consider the composition of the action of $\U'_q(\spa_{2n})$
on $L(\la)$ with the automorphism \eqref{autosic}, where
\ben
\vs_i=\si_i\qquad i=1,3,\dots,2n-1\Fand \vs_i=1,\qquad i=2,4,\dots,2n.
\een
This composition is isomorphic to an irreducible finite-dimensional
highest weight module with the highest weight
$(q^{m_1},\dots,q^{m_n})$. Thus,
without loss of generality,
we may only consider the modules $L(\la)$, where the highest weight
has the form $\la=(q^{m_1},\dots,q^{m_n})$
for some positive integers $m_i$.

Now observe that
for any
$k=1,2,\dots,n-1$,
if we restrict the range of indices of the generators of
$\U'_q(\spa_{2n})$ to the subset $\{2k-1,2k,2k+1,2k+2\}$ then
the corresponding elements will generate a subalgebra
of $\U'_q(\spa_{2n})$ isomorphic to $\U'_q(\spa_{4})$.
The cyclic span of the highest vector $\xi$ of $L(\la)$
with respect to this subalgebra is a highest weight module
with the highest weight $(q^{m_k},q^{m_{k+1}})$.
The irreducible quotient of this module is finite-dimensional.
Hence, considering irreducible highest modules over
$\U'_q(\spa_{4})$ we can get necessary conditions on the $m_i$.

The generators of $\U'_q(\spa_{4})$ are the nonzero entries of the matrix
\ben
S=\left(\begin{matrix} s_{11}&s_{12}&0&0\\
                     s_{21}&s_{22}&0&0\\
                     s_{31}&s_{32}&s_{33}&s_{34}\\
                     s_{41}&s_{42}&s_{43}&s_{44}
\end{matrix}\right)
\een
together with the elements $s^{-1}_{12}$ and $s^{-1}_{34}$.
The highest vector $\xi$ of $L(\la)$ with
the highest weight $\la=(q^{m_1},q^{m_{2}})$
is annihilated
by $s_{11}, s_{31}, s_{32}, s_{33}$ and we have
\ben
s_{12}\ts \xi=q^{m_1}\ts \xi,\qquad
s_{34}\ts \xi=q^{m_2}\ts \xi,
\een
where $m_1$ and $m_2$ are positive integers.
Consider the subspace $L^0$ of $L(\la)$ defined by
\ben
L^0=\{v\in L(\la)\ |\ s_{11}\ts v=s_{31}\ts v=s_{33}\ts v=0\}.
\een
Note that $\xi\in L^0$ and so $L^0\ne 0$.

\ble\label{lem:stable}
The subspace $L^0$ is stable under the action of each of the operators
$s^{}_{32},s^{}_{41},s^{}_{12},s^{}_{34},s^{-1}_{12},s^{-1}_{34}$.
Moreover, these operators
on $L^0$ satisfy the relation
\ben
s^{}_{32}\ts s^{}_{41}-s^{}_{41}\ts s^{}_{32}=(q^2-1)
(s_{12}^{-1}\ts s^{}_{34}-s^{}_{12}\ts s_{34}^{-1}).
\een
\ele

\bpf
The first statement is immediate from \eqref{soopo} for the elements
$s^{}_{12},s^{}_{34}$ and their inverses. The defining
relations \eqref{drabss} imply
\ben
s_{33}\ts s_{32}=s_{32}\ts s_{33},
\qquad
s_{11}\ts s_{32}=s_{32}\ts s_{11}+(q^{-2}-1)\ts s_{12}\ts s_{31},
\een
and
\ben
s_{31}\ts s_{32}=q^{-1}\ts s_{32}\ts s_{31}+(q-\qin)
(q^{-1}\ts s_{21}\ts s_{33}-s_{12}\ts s_{33}),
\een
which implies the statement for the operator $s_{32}$.
Furthermore, we have
\ben
s_{11}\ts s_{41}=s_{41}\ts s_{11},
\qquad
s_{33}\ts s_{41}=s_{41}\ts s_{33}+(q^{-1}-q)\ts s_{34}\ts s_{31},
\een
and
\ben
s_{31}\ts s_{41}=q^{-1}\ts s_{41}\ts s_{31}+(q-\qin)
(q^{-1}\ts s_{43}\ts s_{11}-s_{34}\ts s_{11}),
\een
completing the proof of the first statement.

Now, by \eqref{drabss},
\ben
s_{32}\ts s_{41}=s_{41}\ts s_{32}+(q-\qin)(s_{12}\ts s_{43}-s_{34}\ts s_{21}).
\een
However, \eqref{qdetrel} gives
\ben
s^{}_{21}=q^{-2}\ts s^{-1}_{12}\ts s^{}_{22}\ts s^{}_{11}-q\ts s^{-1}_{12}
\een
so that $s^{}_{21}$ coincides with $-q\ts s^{-1}_{12}$,
as an operator on $L^0$. Similarly,
$s^{}_{43}$ coincides with the operator $-q\ts s^{-1}_{34}$
thus yielding the desired relation.
\epf

Lemma~\ref{lem:stable} implies that $L^0$
is a representation of the quantized enveloping algebra $\U_q(\sll_2)$.
Indeed, the action is defined by setting
\ben
e\mapsto \frac{s_{32}}{q-\qin},\qquad f\mapsto \frac{s_{41}}{q\tss(q-\qin)},
\qquad k\mapsto s^{-1}_{12}\ts s^{}_{34},
\een
where $e,f,k,k^{-1}$ are the standard generators of $\U_q(\sll_2)$
satisfying
\ben
k\tss e=q^2\tss e\tss k, \qquad k\tss f=q^{-2}\tss f\tss k, \qquad
e\tss f-f\tss e=\frac{k-k^{-1}}{q-\qin}.
\een
Since
\ben
e\ts \xi=0\Fand k\ts \xi=q^{-m_1+m_2}\ts \xi,
\een
the cyclic span of $\xi$ with respect to $\U_q(\sll_2)$
is a highest weight module. Since this module is finite-dimensional,
we must have $m_2-m_1\geqslant 0$ due to the classification theorem
for the finite-dimensional irreducible representations of
the algebra $\U_q(\sll_2)$;
see e.g. \cite[Chapter~10]{cp:gq}.
Thus, we have proved the following.

\bpr\label{prop:necco}
Suppose that $\la=(q^{m_1},\dots,q^{m_n})$, where $m_1,\dots,m_n$
are positive integers. If the representation $L(\la)$
of the algebra $\U'_q(\spa_{2n})$ is finite-dimensional
then
$m_1\leqslant m_2\leqslant\dots\leqslant m_n$.
\qed
\epr

Out aim now is to show that these conditions are also sufficient
for the representation $L(\la)$ to be finite-dimensional.
We shall regard $\U'_q(\spa_{2n})$ as a subalgebra of $\U_q(\gl_{2n})$
and use a version of the Gelfand--Tsetlin basis for
representations of $\U_q(\gl_N)$; see \cite{j:qr}.
Let $\nu=(\nu_1,\dots,\nu_N)$
be an $N$-tuple of integers such that $\nu_1\geqslant\dots\geqslant\nu_N$.
The corresponding finite-dimensional irreducible representation
$V(\nu)$ of $\U_q(\gl_N)$
is generated by a nonzero vector $\ze$ such that
\begin{alignat}{2}
\bar t_{ij}\ts\ze&=0 \quad &&\text{for} \quad
1\leqslant i<j\leqslant N,
\non\\
t_{ii}\ts\ze&=q^{\tss\nu_i}\tss\ze \quad
&&\text{for} \quad 1\leqslant i\leqslant N.
\non
\end{alignat}
For any integer $m$ set
\ben
[m]=\frac{q^m-q^{-m}}{q-q^{-1}}.
\een
Define the Gelfand--Tsetlin pattern $\Om$ (associated with $\nu$)
as an array of integer
row vectors of the form
\begin{align}
&\qquad\nu^{}_{N1}\qquad\nu^{}_{N2}
\qquad\qquad\cdots\qquad\qquad\nu^{}_{NN}\non\\
&\qquad\qquad\nu^{}_{N-1,1}\qquad\ \ \cdots\ \
\ \ \qquad\nu^{}_{N-1,N-1}\non\\
&\quad\qquad\qquad\cdots\qquad\cdots\qquad\cdots\non\\
&\quad\qquad\qquad\qquad\nu^{}_{21}\qquad\nu^{}_{22}\non\\
&\quad\qquad\qquad\qquad\qquad\nu^{}_{11}  \non
\end{align}

\bigskip
\noindent
where $\nu^{}_{Ni}=\nu_i$ for $i=1,\dots,N$, so that
the top row coincides with $\nu$, and the following
conditions hold
\ben
\nu_{k+1,i}\geqslant\nu_{ki}\geqslant\nu_{k+1,i+1}
\een
for $1\leqslant i\leqslant k\leqslant N-1$.
There exists a basis $\{\ze^{}_{\tss\Om}\}$ of $V(\nu)$ parameterized
by the patterns $\tss\Om$ such that the action of the generators of
$\U_q(\gl_N)$ is given by
\begin{align}
t_k\ts\ze^{}_{\tss\Om}&=q^{w_k}\ts\ze^{}_{\tss\Om},\qquad w_k=
\sum_{i=1}^k\nu_{ki}-\sum_{i=1}^{k-1}\nu_{k-1,i},
\non
\\
e_k\ts\ze^{}_{\tss\Om}&=-\sum_{i=1}^k\frac{[\ts l_{k+1,1}-l_{ki}]
\cdots [\ts l_{k+1,k+1}-l_{ki}]}
{[\ts l_{k1}-l_{ki}]\cdots\wedge\cdots [\ts l_{kk}-l_{ki}]}
\ts \ze^{}_{\tss\Om+\delta_{ki}},
\non\\
f_k\ts\ze^{}_{\tss\Om}&=\sum_{i=1}^k\frac{[\ts l_{k-1,1}-l_{ki}]
\cdots [\ts l_{k-1,k-1}-l_{ki}]}
{[\ts l_{k1}-l_{ki}]\cdots\wedge\cdots [\ts l_{kk}-l_{ki}]}
\ts \ze^{}_{\tss\Om-\delta_{ki}},
\non
\end{align}
where $l^{}_{ki}=\nu_{ki}-i+1$ and the symbol $\wedge$
indicates that the zero factor in the denominator
is skipped. The array $\Om\pm\delta_{ki}$ is obtained from
$\Om$ by replacing $\nu_{ki}$ with $\nu_{ki}\pm\delta_{ki}$.
The vector $\ze^{}_{\tss\Om}$ is considered to be zero if
the array $\Om$ is not a pattern.

Now let $m_1,\dots,m_n$ be positive integers satisfying
$m_1\leqslant m_2\leqslant\dots\leqslant m_n$.
Consider the representation $V(\nu)$ of $\U_q(\gl_{2n})$
where the highest weight $\nu$ is defined by
\ben
\nu=(r_n,\dots,r_1,0,\dots,0),\qquad r_i=m_i-1.
\een
Introduce the pattern $\Om^{\tss 0}$ associated with $\nu$
which is given by
\begin{align}
&\qquad r_n\qquad r_{n-1}\qquad\cdots\qquad r_1
\qquad 0\qquad\cdots\qquad 0\qquad\quad 0\qquad\quad 0
\non\\
&\qquad\qquad r_n\qquad r_{n-1}\qquad\cdots\qquad r_1
\qquad 0\qquad\cdots\qquad 0\qquad\quad 0
\non\\
&\qquad\qquad\quad\  r_{n-1}\qquad\cdots\qquad r_1
\qquad 0\qquad\cdots\qquad 0\qquad\quad 0
\non\\
&\qquad\qquad\qquad\quad\  r_{n-1}\qquad\cdots\qquad r_1
\qquad 0\qquad\cdots\qquad 0
\non\\
&\qquad\qquad\qquad\qquad\qquad
\cdots\qquad\cdots\qquad\cdots\qquad\cdots\non\\
&\qquad\qquad\qquad\qquad\qquad\qquad
\cdots\qquad\cdots\qquad\cdots\non\\
&\quad\qquad\qquad\qquad\quad\qquad\qquad\qquad r_1\qquad\ 0\non\\
&\quad\qquad\qquad\qquad\quad\qquad\qquad\qquad\qquad r_1 \non
\end{align}

\noindent
For each $k=1,2,\dots,n$ the row $2k-1$ from the bottom
is $(r_k,r_{k-1},\dots,r_1,0,\dots,0)$ with $k-1$ zeros, while
the row $2k$ from the bottom
is $(r_k,r_{k-1},\dots,r_1,0,\dots,0)$ with $k$ zeros.
By the above formulas for the action of the generators of
$\U_q(\gl_{2n})$ in the Gelfand--Tsetlin basis,
for any $i=1,3,\dots,2n-1$ we have the relations
\ben
\bar t_{i,i+1}\ts \ze^{}_{\tss\Om^{\tss 0}}=0\Fand
t_{i,i-1}\ts \ze^{}_{\tss\Om^{\tss 0}}=0,
\een
where the value $i=1$ is excluded for the latter.
The defining relations of $\U_q(\gl_{2n})$ (see Section~\ref{sec:dp})
imply that the vector $\ze^{}_{\tss\Om^{\tss 0}}$
is also annihilated by all generators $t_{jk}$
with $j-k\geqslant 2$ and odd $j$, as well as by
$\bar t_{kj}$ with $j-k\geqslant 2$ and even $j$.

Consider the restriction of the $\U_q(\gl_{2n})$-module $V(\nu)$
to the subalgebra $\U'_q(\spa_{2n})$.
Using \eqref{defrelg2}
and \eqref{sijtij} we then derive that
\ben
s_{ij}\ts \ze^{}_{\tss\Om^{\tss 0}}=0\qquad\text{for} \quad i=1,3,\dots,2n-1,\quad
j=1,2,\dots,i
\een
and
\ben
s^{}_{2k-1,2k}\ts \ze^{}_{\tss\Om^{\tss 0}}=q^{r_k+1}\ts \ze^{}_{\tss\Om^{\tss 0}}
=q^{m_k}\ts \ze^{}_{\tss\Om^{\tss 0}}
\qquad\text{for} \quad k=1,2,\dots,n.
\een
Hence, the cyclic span of the vector $\ze^{}_{\tss\Om^{\tss 0}}$
with respect to the subalgebra $\U'_q(\spa_{2n})$
is a highest weight module with the highest weight
$\la=(q^{m_1},\dots,q^{m_n})$. Since the representation
$V(\nu)$ is finite-dimensional, we conclude
that the representation $L(\la)$ of $\U'_q(\spa_{2n})$
is also finite-dimensional.

Thus, combining this argument
with Propositions~\ref{prop:hwpro} and \ref{prop:necco},
we get the following theorem.

\bth\label{thm:classif}
Every finite-dimensional irreducible representation
of the algebra $\U'_q(\spa_{2n})$ is isomorphic
to a highest weight representation $L(\la)$ where
the highest weight $\la$ is an $n$-tuple of the form
\ben
\la=(\si_1\ts q^{m_1},\dots,\si_n\ts q^{m_n}),
\een
with positive integers $m_i$ satisfying
$m_1\leqslant m_2\leqslant\dots\leqslant m_n$ and each
$\si_i$ is $1$ or $-1$. In particular,
the isomorphism classes of finite-dimensional irreducible representations
of $\U'_q(\spa_{2n})$ are parameterized by such $n$-tuples.
\qed
\eth

It looks plausible that the structure of $L(\la)$ is very much similar
to that of the representation of the Lie algebra $\spa_{2n}$
with the highest weight $(r_n,\dots,r_1)$, where $r_i=m_i-1$.
In particular, these representations should have the same dimensions
and characters. This can be proved
for the case where $\U'_q(\spa_{2n})$ is regarded as an algebra over
$\CC(q)$ by following the arguments of \cite[Section~10.1]{cp:gq}.
Indeed, recall the $\A$-subalgebra $\U'_{\A}$ of $\U'_q(\spa_{2n})$
introduced in Section~\ref{sec:dp}.
Let $\xi$ denote the highest vector of the $\U'_q(\spa_{2n})$-module
$L(\la)$ with $\la=(q^{m_1},\dots,q^{m_n})$, where
the positive integers $m_i$ satisfy
$m_1\leqslant m_2\leqslant\dots\leqslant m_n$.
Set
\ben
L(\la)_{\A}=\U'_{\A}\ts\xi.
\een
Then $L(\la)_{\A}$ is a $\U'_{\A}$-submodule of $L(\la)$
such that
\ben
L(\la)_{\A}\ot^{}_{\A}\ts \CC(q)\cong L(\la)
\een
in an isomorphism of vector spaces over $\CC(q)$.
Moreover, $L(\la)_{\A}$ is the direct sum of its intersections
with the weight spaces of $L(\la)$,
and each intersection
is a free $\A$-module; cf. \cite[Proposition~10.1.4]{cp:gq}.
Now set
\ben
\overline{L}(\la)=L(\la)_{\A}\ot^{}_{\A}\ts\CC,
\een
where the $\A$-action on $\CC$ is defined by the evaluation at $q=1$.
Due to the specialization isomorphism \eqref{speci},
$\overline{L}(\la)$ is a module over the Lie algebra
$\spa_{2n}$. By the specialization formulas of Section~\ref{sec:dp},
we find
\ben
F_{2k-1,2k}\ts\bar\xi=(m_k-1)\ts\bar\xi,\qquad k=1,2,\dots,n,
\een
where $\bar\xi$ denotes the image of $\xi$ in $\overline{L}(\la)$.
Moreover,
\ben
F_{ij}\ts\bar\xi=0,
\qquad\text{for} \quad i=1,3,\dots,2n-1,\quad
j=1,2,\dots,i.
\een
So, taking the weights with respect to
the basis $(F_{2n-1,2n},\dots,F_{12})$ of the Cartan subalgebra
of $\spa_{2n}$, we conclude that $\overline{L}(\la)$
is a highest weight module over $\spa_{2n}$.
Since $\dim\overline{L}(\la)<\infty$, the module $\overline{L}(\la)$
must be irreducible. Thus we come to the following result.

\bth\label{thm:scecia}
The $\spa_{2n}$-module $\overline{L}(\la)$
is isomorphic to the finite-dimensional irreducible
module with the highest weight $(r_n,\dots,r_1)$, $r_i=m_i-1$.
In particular, the character of $\U'_q(\spa_{2n})$-module
$L(\la)$ is given by the Weyl formula
and its dimension over $\CC(q)$ is the same as
the dimension of $\overline{L}(\la)$ over $\CC$.
\qed
\eth


\begin{thebibliography}{99}

\bibitem{acdfr:sb}
{D. Arnaudon, N. Cramp\'e, A. Doikou, L. Frappat, E. Ragoucy},
{\it Spectrum and Bethe ansatz equations
for the $U_q(gl(N))$ closed and
open spin chains in any representation},
preprint {\tt math-ph/0512037}.

\bibitem{cp:gq}
{V. Chari and A. Pressley},
{\it A guide to quantum groups},
Cambridge University Press, 1994.

\bibitem{d:ha}
{V. G. Drinfeld},
{\it Hopf algebras and the
quantum Yang--Baxter equation},
{Soviet Math. Dokl.} {\bf
32} (1985), 254--258.

\bibitem{gk:qd}
{A. M. Gavrilik and A. U. Klimyk},
{\it $q$-deformed orthogonal and pseudo-orthogonal algebras and their
representations}, Lett. Math. Phys. {\bf 21} (1991), 215--220.

\bibitem{gi:ce}
{A. M. Gavrilik and N. Z. Iorgov},
{\it On Casimir elements of $q$-algebras
$U'_q({\rm so}_n)$ and their
eigenvalues in representations}, in
`Symmetry in nonlinear mathematical physics',
Proc. Inst. Mat. Ukr. Nat. Acad. Sci.
{\bf 30},
Kyiv, 1999, pp. 310--314.

\bibitem{gik:nd}
{A. M. Gavrilik, N. Z. Iorgov and A. U. Klimyk},
{\it Nonstandard deformation $U'_q({\rm so}_n)$: the embedding
$U'_q({\rm so}_n)\subset U_q({\rm sl_n})$ and representations}, in
`Symmetries in science', X (Bregenz, 1997).
Plenum, New-York, 1998, pp. 121--133.

\bibitem{hkp:ce}
{M. Havl\'\i\v cek, A. U. Klimyk and S. Po\v sta},
{\it Central elements of the algebras $U'\sb q({\rm so}\sb m)$
and $U\sb q({\rm iso}\sb m)$},
Czechoslovak J. Phys. {\bf 50} (2000), 79--84.

\bibitem{h:il}
{J. E. Humphreys},
{\it Introduction to Lie algebras and Representation Theory},
Springer, New York, 1972.

\bibitem{ik:nd}
{N. Z. Iorgov and A. U. Klimyk},
{\it The nonstandard
deformation $U'_q({\rm so}_n)$ for $q$ a root of unity},
Methods of Funct. Anal. Topology {\bf 6}, (2000), 15--29.

\bibitem{ik:ct}
{N. Z. Iorgov and A. U. Klimyk},
{\it Classification theorem on irreducible representations
of the $q$-deformed algebra $\U'_q(\mathfrak{so}_n)$},
Int. J. Math. Sci.  2005,  no. 2, 225--262.

\bibitem{j:qd}
{M. Jimbo},
{\it A $q$-difference analogue of $\U(\g)$ and the Yang--Baxter equation},
Lett. Math. Phys. {\bf 10} (1985), 63--69.

\bibitem{j:qu}
{M. Jimbo},
{\it A $q$-analogue of $U_q(\gl(N+1))$, Hecke algebra and
the Yang--Baxter equation},
Lett. Math. Phys. {\bf 11} (1986), 247--252.

\bibitem{j:qr}
{M. Jimbo},
{\it Quantum $R$-matrix for the generalized Toda system},
{Comm. Math. Phys.}
{\bf 102} (1986),
537--547.

\bibitem{ks:qg}
{A. Klimyk and K. Schm\"{u}dgen},
{\it Quantum qroups and their
representations}, Springer, Berlin, 1997.

\bibitem{l:sp}
{G. Letzter},
{\it Symmetric pairs for quantized enveloping algebras},
{J. Algebra} {\bf 220} (1999), 729--767.

\bibitem{l:cs}
{G. Letzter},
{\it Coideal subalgebras and quantum symmetric pairs},
in `New directions in Hopf algebras',
Math. Sci. Res. Inst. Publ. {\bf 43},
Cambridge Univ. Press, Cambridge, 2002, pp. 117--165.

\bibitem{m:nq}
{A. I. Molev},
{\it A new quantum analog of the Brauer algebra}, {Czech.
J. Phys.} {\bf 53} (2003), 1073--1078.

\bibitem{mrs:cs}
{A. Molev, E. Ragoucy and P. Sorba},
{\it Coideal subalgebras in
quantum affine algebras}, Rev. Math. Phys.
{\bf 15} (2003), 789--822.

\bibitem{n:ms}
{M. Noumi},
{\it Macdonald's symmetric polynomials as zonal spherical
functions on quantum homogeneous spaces},
Adv. Math. {\bf 123} (1996), 16--77.

\bibitem{nuw:dp}
{M. Noumi, T. Umeda and M. Wakayama},
{\it Dual pairs, spherical harmonics and a Capelli identity
in quantum group theory},
{Compos. Math.}
{\bf 104} (1996),
227--277.

\bibitem{rtf:ql}
{N. Yu. Reshetikhin, L. A. Takhtajan and L. D. Faddeev},
{\it Quantization of Lie Groups and Lie algebras}, {Leningrad Math. J.}
{\bf 1}
(1990),
193--225.

\bibitem{r:ap}
{M. Rosso},
{\it An analogue of P.B.W. theorem and the universal $R$-matrix
for $U_h sl(N+1)$},
Comm. Math. Phys. {\bf 124} (1989), 307--318.


\end{thebibliography}
\end{document}